\documentclass{article}

\usepackage{epsfig}
\usepackage{latexsym}

\newtheorem{theorem}{Theorem}[section]
\newtheorem{corollary}[theorem]{Corollary}

\newtheorem{question}[theorem]{Question}
\newtheorem{lemma}[theorem]{Lemma}

\newenvironment{proof}{{\bf Proof.}}{\hfill\rule{2mm}{2mm}}

\def\mod{ {\rm mod \ }}

\newtheorem{prelem}{{\bf Theorem}}

\newtheorem{prelemm}{{\bf Lemma}}

\def\Ex {{\bf E}}

\title{\bf Random cubic graphs are not homomorphic to the cycle of size $7$.}
\author{
{\bf  Hamed Hatami} \\
{\small\it Department of Computer Science}\\
{\small University of Toronto} \\
{\small e-mail: hamed@cs.toronto.edu}}
\date{}
\begin{document}
\maketitle
\begin{abstract}
We prove that  a random cubic graph almost surely is not
homomorphic to a cycle of size $7$. This implies that there exist
cubic graphs of arbitrarily high girth with no homomorphisms to
the cycle of size $7$.
\end{abstract}
{{\sc Keywords:} homomorphism, cubic graph, random cubic graph,
circular chromatic number.

\section{Introduction}        
For a graph $G$, we denote its vertex set by $V(G)$. Suppose $G$
and $H$ are graphs. A {\sf homomorphism} from $G$ to $H$ is a
mapping $h$ from $V(G)$ to $V(H)$  such that for each edge $xy$ of
$G$, $h(x)h(y)$ is an edge of $H$. We say that $G$ is homomorphic
to $H$, if there exists a homomorphism from $G$ to $H$. A
homomorphism from $G$ to $K_n$ is equivalent to an $n$-coloring of
$G$.  So graph homomorphism is a generalization of coloring. Since
every even cycle is bipartite, it is homomorphic to a single edge.
So for every even $n>0$, a graph $G$ is homomorphic to $C_n$, if
and only if $G$ is bipartite.

Suppose that $G \neq K_4$ is a cubic graph. Then $G$ is
$3$-colorable, and so it is homomorphic to $C_3$. Since a graph
containing $C_i$ for $i$ odd has no homomorphism onto $C_j$ for
any $j>i$, if $G$ contains a triangle, then it is not homomorphic
to $C_5$. The Peterson graph is triangle-free, but it is not
homomorphic to $C_5$. The following question is first asked
in~\cite{nesetril1} (see also~\cite{nesetril}).

\begin{question}
\label{cycle5} Is it true that any cubic graph $G$ with
sufficiently large girth $g$ is homomorphic to $C_5$?
\end{question}

It is shown in~\cite{nesetril} that the answer is negative when
$C_5$ is replaced with $C_{11}$. In~\cite{wormald} Wanless and
Wormald improved this result to $C_9$ by studying the problem on
random cubic graphs. In this note we improve their result to $C_7$
by showing that a random cubic graph almost surely is not
homomorphic to $C_7$, where we say that a property $P(n)$ {\sf
a.s.} (almost surely) holds,  if $\Pr[P(n)]=1-o(1)$. It is easy to
see that there exist nonbipartite cubic graphs of arbitrarily high
girth. So we leave $C_5$ as the only remaining open case.

We consider the following probability space. Choose three random
perfect matchings independently and uniformly on $n$ ($n$ even)
vertices. Let $G$ be the multigraph obtained by taking the union
of these three perfect matchings.  If $G$ is restricted to having
no multiple edges, then any property holds for $G$ a.s., if and
only if it holds a.s. for a random cubic graph
(see~\cite{wormaldsurvey}). The probability that $G$ has no
multiple edges is asymptotic to $e^{-\frac{3}{2}}<1$
(see~\cite{wormald}). Thus we can conclude that if a property
holds for $G$ a.s., then it holds a.s. for a random cubic graph.

Circular chromatic number is a generalization of the chromatic
number. It is not known if there exists a cubic  graph with
arbitrarily large girth whose circular chromatic number is exactly
$3$. Our result shows that there exist cubic graphs with
arbitrarily large girth whose circular chromatic number is more
than $7/3=2.33\ldots$.

In Section~2 we prove that a random cubic graph a.s. is not
homomorphic to the cycle of size $7$. In Section~3 we introduce
the relation between our result and the circular chromatic number
of cubic graphs, and we pose some open problems.

\section{Main result}

We take $C_k$ with vertex set being the congruence classes modulo
$k$, and edges joining $i$ to $i+1$, where in such notation the
integers represent their congruence classes. Let $h$ be a
homomorphism from a graph $G$ to $C_k$. We say that $h$ is {\sf
tight}, if for every vertex $v$ where $h(v) \neq 0$, there exists
a vertex $u$ adjacent to $v$ such that $h(u) = h(v)+1$.

\begin{lemma}
\label{jackpot} Suppose that $G$ is a graph on $n$ vertices which
has a homomorphism to $C_{2k+1}$. Then there is a tight
homomorphism from $G$ to $C_{2k+1}$.
\end{lemma}
\begin{proof}
Begin with a homomorphism $h_0$ from $G$ to $C_{2k+1}$ which maps
every isolated vertex of $G$ to $0$. For every $i\ge 0$, if $h_i$
is not tight, then the homomorphism $h_{i+1}$ is defined
recursively as in the following. Let $v$ be a vertex in $G$ such
that $h_i(v) \neq 0$, and $h_i(u)=h_i(v)-1$ for every vertex $u$
adjacent to $v$. Define $h_{i+1}$ as follows:
$h_{i+1}(v)=h_i(v)-2$, and for every other vertex $w$,
$h_{i+1}(w)=h_i(w)$. Observe that $h_{i+1}$ is a homomorphism, and
$h_i(v)=0$ implies that $h_j(v)=0$ for every $j \ge i$. Since for
every $0 \le x \le 2k$ there exists $0 \le i \le 2k$ such that
$x-2i=0 \ (\mod 2k+1)$, there is some $t$ such that $h_t$ is a
tight homomorphism. In fact $t<2kn$.
\end{proof}

In~\cite{wormald} the expected number of homomorphisms of a random
cubic graph to $C_9$ is studied. They showed that this expected
number is $o(1)$ conditioned on some property which is true,
almost surely. For $C_7$ this value is exponentially high, so
instead of the expected number of homomorphisms, we focus on the
expected number of tight homomorphisms. This allows us to improve
the previous result.

\begin{theorem}
\label{main} A random cubic graph a.s. is not homomorphic to $C_7$.
\end{theorem}
\begin{proof}
Let $M_1$, $M_2$, and $M_3$ be three random perfect matchings
chosen independently and uniformly on $n$ vertices, and $G$ be the
multigraph obtained by taking the union of $M_1$, $M_2$, and
$M_3$. Let $A$ denote the event that the size of the largest
independent set of a graph is less than $0.4554n$.
McKay~\cite{McKay} showed that $A$ holds a.s. for every cubic
graph which implies that $A$ holds a.s. for $G$. Let $I_A$ be the
indicator variable of $A$. We will bound $\Ex[X(G) I_A(G)]$, where
$X(G)$ is the number of tight homomorphisms from $G$ to $C_7$. If
we can show that $\Ex[X(G) I_A(G)]=o(1)$, then since a.s.
$I_A(G)=1$, we can conclude that $X(G)=0$ almost surely. Then
Lemma~\ref{jackpot} implies that a.s. $G$ does not have any
homomorphism to $C_7$.

Let ${\cal T}$ be the set of ordered triples of perfect matchings
on $n$ vertices, and ${\cal H}$ be the set of all mappings from
$n$ vertices to $C_7$. For every $h \in {\cal H}$, let $t(h)$ be
the number of triples $T \in {\cal T}$ such that their
corresponding graph $G$ has property $A$ and $h$ is a tight
homomorphism from $G$ to $C_7$. We have

$$\Ex[X(G) I_A(G)] = \frac{1}{|{\cal T}|} \sum_{h \in {\cal H}} t(h).$$

Consider a triple $T$ which contributes to $t(h)$ for some mapping
$h$. For each $i$, let $n_i$ be the cardinality of $h^{-1}(i)$.
Consider a particular matching $M$ in $T$, and let $m_i$ be the
number of edges between $h^{-1}(i-1)$ and $h^{-1}(i)$ in $M$. Then
the following trivial equation shows that the values of $m_i$ are
equal for all of the three matchings in $T$.
\begin{equation}
\label{eq1} m_i=\frac{n}{2}-(n_{i+1} +n_{i+3} + n_{i+5}).
\end{equation}
Note that since $h^{-1}(i+1) \cup h^{-1}(i+3) \cup h^{-1}(i+5)$ is
an independent set in $G$, we have $n_{i+1} +n_{i+3} + n_{i+5} \le
0.4554n$ which implies that $m_i \ge 0.0446n$.

Let us focus on $h^{-1}(i)$ for some $i>0$. Each perfect matching
in $T$ partitions $h^{-1}(i)$ into two sets of size $m_i$ and
$m_{i+1}$ which are the vertices that are adjacent to
$h^{-1}(i-1)$ and $h^{-1}(i+1)$, respectively in that perfect
matching. Since $h$ is a tight homomorphism, no element of
$h^{-1}(i)$ is in the partition of size $m_i$ in all of the three
perfect matchings. So the number of different possible ways to
partition the sets is
$${n_0 \choose m_0}^3 \prod_{i=1}^{6}  \sum_{j=0}^{m_i} {n_i \choose m_i} {m_i \choose j}
{n_i-m_{i} \choose m_i-j} {n_i-j \choose m_i}.$$

Note that ${a \choose b}$ is defined to be $0$, for $b>a$. When
the partitions are determined, there are $\prod_{i=0}^{6}m_i!^3$
different ways to arrange matching edges on them. Since there are
$\frac{n!}{n_0!n_1!\ldots n_6!}$ different ways to partition the
vertices into groups of size $n_0,n_1,\ldots,n_6$, we have
$$\sum_{h \in {\cal H}} t(h) \le \sum_{n_0+\ldots+n_6=n}\frac{n!}{n_0!\ldots n_6!}
\left(\prod_{i=0}^{6}m_i!^3\right) {n_0 \choose m_0}^3
\prod_{i=1}^{6} \sum_{j=0}^{m_i} {n_i \choose m_i} {m_i \choose j}
{m_{i+1} \choose m_i-j} {n_i-j \choose m_i}.$$
The outer sum has less than $n^6$ summands, and each one of the
six inner sums are taken over at most $n$ terms. So by
substituting $n_i=m_i+m_{i+1}$ and converting sums to maximums we
have:
\begin{equation}
\label{eq2} \Ex[X(G) I_A(G)] \le \frac{n!}{|{\cal T}|} \times
n^{12} \max_{m_0,\ldots,m_6} f(m_0,m_1,\ldots,m_6),
\end{equation}
where
$$f(m_0,m_1,\ldots,m_6)={m_0+m_1 \choose m_0}^2\left(
\prod_{i=0}^{6}m_i!\right) \prod_{i=1}^{6} \max_{0 \le j \le m_i}
{m_i \choose j} {m_{i+1} \choose m_i-j} {m_i+m_{i+1}-j \choose
m_i}.$$
Note that if $j>m_{i+1}$ or $j<m_i-m_{i+1}$, then ${m_i \choose j}
{m_{i+1} \choose m_i-j} {m_i+m_{i+1}-j \choose m_i}=0$. Define
$y_i=m_i/n$ and
$$B_i=[0, y_i] \cap [0, y_{i+1}]\cap [y_{i}-y_{i+1},\infty].$$
Let $g(x)=x^x$, for $x>0$; and $g(0)=1$.
 Then
\begin{equation}
\label{eq3} f(m_0,\ldots,m_6)^{1/n}=\sqrt{\frac{n}{e}}
\left(\frac{g(y_0+y_1)^2}{g(y_0) g(y_1)^2} \prod_{i=1}^{6} \max_{
z \in B_i }
\frac{g(y_i)g(y_{i+1})g(y_i+y_{i+1}-z)}{g(z)g(y_i-z)^2g(y_{i+1}-y_i+z)g(y_{i+1}-z)}+o(1)\right).
\end{equation}

 Note that if $y_i > 2y_{i+1}$ for some
$0<i\le 6$, then $B_i$ is empty and $f(m_0,\ldots,m_6)=0$. By
substituting $|{\cal T}|=\left( \frac{n!}{2^{n/2}(n/2)!}\right)^3$
in Equation~(\ref{eq2}), we have
\begin{equation}
\label{maineq} \Ex[X(G) I_A(G)]^{1/n} \le \max
\frac{g(y_0+y_1)^2}{g(y_0) g(y_1)^2} \prod_{i=1}^{6}\max_{ z \in
B_i }
\frac{g(y_i)g(y_{i+1})g(y_i+y_{i+1}-z)}{g(z)g(y_i-z)^2g(y_{i+1}-y_i+z)g(y_{i+1}-z)}+o(1),
\end{equation}
where the outer maximum is taken with respect to the following
conditions,
\begin{equation}
\label{condition}
\begin{array}{rclc}
\sum_{i=0}^{6} y_i  & = & \frac{1}{2} & \\
y_i & \le & 2y_{i+1} & (1 \le i \le 6) \\
y_i  & \ge & 0.0446  & (0 \le i \le 6)
\end{array}
\end{equation}

 Suppose that for $0 \le i \le 6$, $0 \le a_i \le b_i$
are given, where $b_i-a_i \le 0.1$. We want to examine Inequality~(\ref{maineq}) when $a_i
\le y_i \le b_i$, for every $0 \le i \le 6$. To satisfy
Condition~(\ref{condition}) we have the following restrictions on
$a_i$ and $b_i$.

\begin{equation}
\label{condition2}
\begin{array}{cc}
\sum_{i=0}^{6} a_i   \le  \frac{1}{2} \le  \sum_{i=0}^{6} b_i& \\
a_i  \le  2b_{i+1} & (1 \le i \le 6) \\
b_i   \ge  0.0446  & (0 \le i \le 6) \\
\end{array}
\end{equation}

By considering the
derivatives of $\ln\left(\frac{g(x+y)^2}{g(x)
g(y)^2}\right)$ with respect to $x$ and $y$,
we conclude that
$$\frac{g(y_0+y_1)^2}{g(y_0) g(y_1)^2} \le \frac{g(y_0+b_1)^2}{g(y_0) g(b_1)^2} \le
\max(\frac{g(a_0+b_1)^2}{g(a_0) g(b_1)^2},\frac{g(b_0+b_1)^2}{g(b_0) g(b_1)^2})=f_0(a_0,b_0,b_1).$$
If Condition~(\ref{condition}) is satisfied, then obviously for every
$0 \le i \le 6$, we have $y_i \le 0.5-6\times 0.0446=0.2324$. So $b_{i+1}-a_i+z<\frac{1}{e}=0.367\ldots$, which
implies that for $z \in B_i$,
$$\frac{g(y_i)g(y_{i+1})g(y_i+y_{i+1}-z)}{g(z)g(y_i-z)^2g(y_{i+1}-y_i+z)g(y_{i+1}-z)} \le
\frac{g(b_i)g(b_{i+1})g(b_i+b_{i+1}-z)}{g(z)g(b_i-z)^2g(b_{i+1}-z)}\times
\frac{1}{g(b_{i+1}-a_i+z)}.$$
Since $B_i \subseteq B_i'$, where $B_i'$ is defined as
$$B_i'=[0, b_i] \cap [0,b_{i+1}]\cap [a_{i}-b_{i+1},\infty],$$
we conclude that
\begin{equation}
\label{maxcompare} \max_{z \in B_i}
\frac{g(y_i)g(y_{i+1})g(y_i+y_{i+1}-z)}{g(z)g(y_i-z)^2g(y_{i+1}-y_i+z)g(y_{i+1}-z)}
\le \max_{z \in B_i'}
\frac{g(b_i)g(b_{i+1})g(b_i+b_{i+1}-z)}{g(z)g(b_i-z)^2g(b_{i+1}-a_i+z)g(b_{i+1}-z)}.
\end{equation}

The derivative with respect to $z$ of the natural logarithm of the expression
following $\max$ on the right hand side of Inequality~(\ref{maxcompare}) is
$$\ln\left(\frac{(b_{i+1}-z)(b_i-z)^2}{z(b_i+b_{i+1}-z)(b_{i+1}-a_i+z)}\right).$$
The only critical value of $z\in B'_i$ is
$z_i=(-B-\sqrt{B^2-4AC})/(2A),$ where
$$A=b_{i+1}+b_i-a_i,$$
$$B=b_{i+1}a_i+b_{i}a_i-3b_{i+1}b_i-b_i^2-b_{i+1}^2,$$
$$C=b_i^2b_{i+1}.$$
Since the derivative tends to infinity as $z$ tends from above to
$\max(0,a_i-b_{i+1})$, and it tends to minus infinity as $z$ tends
 from below
to $\min(b_i,b_{i+1})$, the function is maximized at $z_i$.
So finally we have

\begin{eqnarray*}
\lefteqn{\max_{z \in B'_i}\frac{g(y_0+y_1)^2}{g(y_0) g(y_1)^2}
\prod_{i=1}^{6}\max_{ z \in B_i }
\frac{g(y_i)g(y_{i+1})g(y_i+y_{i+1}-z)}{g(z)g(y_i-z)^2g(y_{i+1}-y_i+z)g(y_{i+1}-z)}
\le} \\ & & f_0(a_0,b_0,b_1) \prod_{i=1}^{6}
\frac{g(b_i)g(b_{i+1})g(b_i+b_{i+1}-z_i)}{g(z_i)g(b_i-z_i)^2g(b_{i+1}-a_i+z_i)g(b_{i+1}-z_i)}=h(a_0,\ldots,a_6,b_0,\ldots,b_6).
\end{eqnarray*}

Fix some $0<\epsilon<0.1$, and let $a_i= \lfloor
\frac{y_i}{\epsilon} \rfloor \epsilon$, and $b_i=a_i+\epsilon$.
Since $a_i$ is the product of $\epsilon$ and some integer, there
are finite number of different possible values for each $a_i$. A
computer is used to compute $h(a_0,\ldots,a_6,b_0,\ldots,b_6)$,
for these values when Condition~(\ref{condition2}) is satisfied.
It turned out that for $\epsilon=0.00125$ we have
$h(a_0,\ldots,a_6,b_0,\ldots,b_6)< 0.99$ for all possible values
of $a_i$ and $b_i$ which satisfy Condition~(\ref{condition2}).
This implies that
 $$\Ex[X(G) I_A(G)]^{1/n} < 0.99+o(1),$$
 which is
sufficient to complete the proof.
\end{proof}

\begin{remark}
For $\epsilon=0.00125$, the number of possible values of $a_i$ and
$b_i$ is very large, and  it takes a long time for a computer to
do the required computations. To avoid this problem, we started
with $\epsilon_1=0.01$. Every time that some $a_i$ and $b_i$,
$b_i=a_i+\epsilon_k$ ($0 \le i \le 6$) satisfying
Condition~(\ref{condition2}) and
$h(a_0,\ldots,a_6,b_0,\ldots,b_6)>0.99$ are found, we computed $h$
for all $a_i \le a'_i \le b'_i \le b_i$ ($0 \le i \le 6$) where
$a'_i$ is the product of $\epsilon_{k+1}$ and some integer,
$\epsilon_{k+1}=\frac{\epsilon_k}{2}$, $b'_i=a_i'+\epsilon_{k+1} $
and Condition~(\ref{condition2}) is satisfied. It turned out that
$\epsilon_4=0.00125$ is always sufficient.
\end{remark}

For every integer $g$ there is a constant $\epsilon$ such that the
probability that the girth of a random cubic graph is greater than
$g$ is more than $\epsilon$ (see for
example~\cite{wormaldsurvey}, Section 2.3). So we have the
following corollary.

\begin{corollary}
For every integer $g$, there is a cubic graph $G$ with girth at
least $g$ which is not homomorphic to the cycle of size $7$.
\end{corollary}

\section{Applications to the circular chromatic number}

For a pair of integers $p$ and $q$ such that $p \ge 2q$, let
$K_{p/q}$ be the graph that has vertices $\{0,1,\ldots,p-1\}$ and
in which $x$ and $y$ are adjacent if and only if $q \le |x-y| \le
p-q$. For a graph $G$, the circular chromatic number $\chi_c(G)$
is defined as the infimum of $p/q$ where there is a homomorphism
from $G$ to $K_{p/q}$. It is known that for every graph $G$ the
infimum in the definition is always attained and $\chi(G)-1 <
\chi_c(G) \le \chi(G)$ (see for example~\cite{zhusurvey}). Since
$K_{7/3}$ is a cycle of size $7$, Theorem~\ref{main} shows that
the circular chromatic number of a random cubic graph is a.s.
greater than $7/3$ which implies the following corollary.
\begin{corollary} There exist cubic graphs of arbitrary large girth whose
circular chromatic number is greater than $7/3$.
\end{corollary}
We ask the following question.
\begin{question}
Determine the supremum of $r^*$ such that there exist cubic graphs
of arbitrary large girth whose circular chromatic number is at
least $r^*$.
\end{question}

We know that $7/3 \le r^* \le 3$. Since the circular chromatic
number of $C_5$ is $5/2$, the answer of Question~\ref{cycle5}
determines whether $r^* \ge 5/2$ or $r^* \le 5/2$. If the answer
to Question~\ref{cycle5} is affirmative, then it would be easier
to first answer the following question.

\begin{question}
Is it true that the circular chromatic number of any cubic graph
$G$ with sufficiently large girth  is less than $3$?
\end{question}

\section*{Acknowledgment}
The author wishes to thank Michael Molloy for his valuable
discussions leading towards the result of this paper and his
valuable comments on the draft version, and Xuding Zhu for
introducing the problem.
\bibliographystyle{siam}
\bibliography{hom}

\end{document}